
\magnification 1200
\def\qed{\hskip .6em \raise1.8pt\hbox{\vrule height4pt width6pt depth2pt}}
\def\implies{\hbox{$\Rightarrow$}}
\font\Bbb=msbm10

\def\N{\hbox{\Bbb N}}
\baselineskip 18pt

\centerline
{{\bf When each continuous operator is regular, II} \footnote*
{This work was started whilst the second named author was visiting Indiana
University-Purdue University at Indianapolis in the summer of 1993 and
finished during the visit of the first named author to the Queen's
University of Belfast in the spring of 1995, under  the auspices of
a NATO Collaborative Research Grant CRG--890909.}}
\bigskip

\centerline {\bf Y. A. Abramovich and A. W. Wickstead }
\bigskip

{\narrower\narrower {\bf Abstract.} The following theorem is essentially due to
L.~Kantorovich and B. Vulikh and it describes  one of the most important
classes of Banach lattices between which each continuous operator is regular.
{\bf Theorem 1.1.} {\sl Let $E$ be an arbitrary L-space and $F$ be an arbitrary
Banach lattice with Levi norm. Then ${\cal L}(E,F)={\cal L}^r(E,F),\ (\star) $
that is, every continuous operator from $E$ to $F$ is regular.}

In spite of the importance of this theorem
it has not yet been determined  to what extent  the Levi condition is
essential
for the validity of equality $(\star)$. Our main aim in  this work is to prove a
converse  to this theorem  by showing that for a Dedekind complete  $F$
the Levi condition is necessary for the validity of $(\star)$.

As a sample of other results we mention the following. {\bf Theorem~3.6.} {\sl
For a Banach lattice $F$  the following are equivalent:
{\rm (a)} $F$ is Dedekind complete;
{\rm (b)} For all Banach lattices $E$, the space ${\cal L}^r(E,F)$ is a Dedekind
complete  vector lattice;
{\rm (c)} For all L-spaces $E$, the space ${\cal L}^r(E,F)$ is a vector
lattice.}

}

\bigskip

\noindent{\bf 1. Introduction.}  As the title of this work indicates we will
be concerned here with the study of Banach lattices $E$ and $F$ for which
the space of all continuous operators, ${\cal L}(E,F)$, coincides with the space
of all regular operators, ${\cal L}^r (E,F)$. Recall that a (linear)
operator is said to be regular if it can be split into the difference of two
positive operators. The following theorem, which is essentially due to
L. Kantorovich and B. Vulikh  [KV], describes  one of the most important
classes of
Banach lattices between which each continuous operator is regular.

\medskip
\noindent{\bf Theorem 1.1.}
{\sl Let $E$ be an arbitrary L-space and $F$ be an arbitrary
Banach lattice with Levi norm. Then
             $${\cal L}(E,F)={\cal L}^r(E,F), \eqno(\star) $$
that is, every continuous operator from $E$ to $F$ is regular.}

\medskip
To be precise, it was assumed in [KV] that $F$ was a KB-space,
and it was noticed in [S] that the original proof could be easily carried over
from a KB-space $F$ to an arbitrary   Banach lattice with a Levi norm.
For a KB-space $F$ the proofs can be  found in [AB, Theorem~15.3] and
[V, Theorem~8.7.2]. Under the assumption, somewhat stronger than Levi
property, that $F$ is positively complemented in $F^{**}$ a proof of
Theorem~1.1 is presented in [MN, Theorem 1.5.11].

It is somewhat  surprising that  in spite of the importance of Theorem~{1.1}
it has not yet been determined  to what extent  the Levi condition is
essential
for the validity of equality $(\star)$. Our main aim in  this work is to prove a
converse  to Theorem~{1.1} by showing that for a Dedekind complete  $F$
the Levi condition is necessary for the validity of $(\star)$.

\medskip
\noindent{\bf Definition 1.2.} {\sl A norm on a Banach lattice $E$
is said to be a  {\rm Levi norm}, if every norm-bounded
upward directed set of positive elements has a supremum.
If the previous property holds only for sequences, then
we say that the norm is {\rm sequentially Levi.}}

Obviously, each  Banach lattice with a  Levi norm must be Dedekind
complete, and  each Banach  lattice with a  sequentially Levi norm must
be Dedekind $\sigma$-complete. It is worth noticing that the above
definition is, in fact,  of an order-topological nature as it describes
relationships between the topology and the order, rather than
the properties of a particular (lattice) norm.  These properties
appear in the literature under many different names. It was
D.~Fremlin [F] who was the first to use Levi's name in connection
with this property. H.~Nakano [N, pages~129-130] used the term  monotone
complete norm for a sequentially Levi norm and universally monotone complete
norm for Levi.
A. Zaanen [Z] considers sequentially Levi norms in two places under
different names. First on page 305 he refers to them, like Nakano,
as monotone complete norms, and then on page 421 as norms with the weak
Fatou property for monotone sequences. The term weak Fatou property
for directed sets (page 390) is used by Zaanen for what we refer to
as a Levi norm.
P.~Meyer-Nieberg [MN, page~96]
uses  the term monotonically complete. Finally, the Soviet school on Banach
lattices  used symbols (B) and (B$'$) to denote  sequentially Levi
and Levi properties respectively.

Our converse to Theorem~{1.1}, which we mentioned above, is somewhat partial
since we assume  $F$ to be Dedekind complete. On the other hand, it is the
best one can get since, as we will see below in Remark~3.2 and the comments
after Theorem~3.5, there are  non-Dedekind complete Banach lattices (hence,
without  a Levi norm)  which, nevertheless, satisfy the equality
$(\star)$.

This paper can be considered as a sequel to [A3]. As in [A3] we
adhere in this work to an isomorphic point of view, i.e., we do
not  distinguish between equivalent norms.
% A single isometric result will be presented in Section~4.
We use the standard
terminology regarding Banach lattices and operators on them.
Any notation or definition not mentioned explicitly in the
text can be found in [AB], [V] or [A3].

\bigskip
\noindent{\bf 2. Some Banach lattice preliminaries.} In this section
we  present two new results of the so-called lateral analysis
which will be needed later on. Lateral analysis is a  convenient and
powerful method  of investigating Banach and vector lattices. The
essence of this method can be  described roughly as follows: instead of
arbitrary nets or sequences (to be considered in a property or a
definition)  we try to deal  with those of a much simpler structure by
considering  only the nets or sequences with mutually  disjoint terms.  We
refer to  [A2], [AB] and [MN] where this approach is used systematically.
For an illustration we present  two examples
 the former of which will be used later on.

The theorem of Veksler and Gejler [VG],
characterizes Dedekind completeness of vector lattices
by  stating  that a  uniformly complete  vector
lattice is Dedekind complete if and only if every order
bounded set of pairwise disjoint positive elements
has a supremum. Many other completeness  properties of vector
or  Banach lattices  have been also characterized in the
framework of lateral  analysis. Meyer-Nieberg [MN1] and Fremlin
[F, page~56] have shown that a Banach lattice has an order
continuous norm if and only if every order bounded  sequence of pairwise
disjoint  elements converges to zero in norm. (See also [AB, Theorem 12.13]
or [MN, Theorem 2.4.2] for alternative proofs of this theorem.)

Recall that a vector lattice  is said to be
{\it universally complete\/}
if it is Dedekind complete and has the property that every set
of pairwise disjoint positive elements has a supremum.
It is well known [V, Chapter V] that every Dedekind complete
vector lattice $E$ has a universal completion, $\hat E$,
which, by definition,  is a universally complete vector lattice
containing $E$ as an order dense ideal.

\medskip
\noindent{\bf Definition 2.1.} {\sl If $E$ is a vector lattice  then an
upward  directed set
$A\subseteq E_+$ is called {\rm laterally
increasing}  if  for each  $a,b\in A$ with $a\ge b$ we have
$(a-b)\wedge b=0$.}

There is  an important difference between laterally increasing
sequences and nets which we would like to point out.
If $(x_n)$ is  a laterally increasing sequence, then one can
easily produce a sequence $(u_n)$ with pairwise disjoint elements
such that $x_n=u_1+\ldots+u_n$ (take simply $u_1=x_1$ and
$u_n=x_{n+1}-x_{n}$ for $n\ge 2 $).  For laterally increasing nets,
however, there is no convenient substitute  for the previous
representation, and this makes  working with nets more complicated.
Our next proposition and theorem deal with this problem.

\medskip
\noindent{\bf Proposition 2.2.} {\sl A Dedekind complete vector lattice
$E$ is universally complete if and only if every laterally increasing
subset of
$E_+$ has a supremum.}

\smallskip
\noindent{\bf Proof.} If every laterally increasing subset of
$E_+$ has a supremum and $A\subset E_+$ is a  given
pairwise disjoint set, then the set of all finite sums of elements from
$A$ is laterally increasing so has a supremum. That supremum is clearly
also the supremum of $A$ itself, so that $E$ is indeed universally
complete.

Now suppose that $E$ is universally complete. We know that $E$
is isomorphic to a space $C_\infty(Q)$ for some Stonean space $Q$ [V,
Chapter V]. Suppose  that $A$ is a
laterally increasing subset of $C_\infty(Q)$. Define a function $y$ on
$Q_0=\bigcup_{a\in A}\{s\in Q: a(s)>0\}$ by
$y(s)=a(s)$ if $a(s)>0$. In order to show that this definition is
unambiguous, it suffices to consider $a\ge b$ with $b(s)>0$ and show
that we obtain the same value for $y(s)$ using either $a$ or $b$, for
then if we have any $b,c\in A$ we need only take $a\ge b,c$ to see that
$b$ gives the same value as $a$, which in turn gives the same value as
$c$. But if $a\ge b$ then $(a-b)\wedge b=0$, and in particular
$(a-b)(s)\wedge b(s)=0$. As $b(s)>0$ this means that $(a-b)(s)=0$, i.e.
$a(s)=b(s)$ and the definition is therefore unambiguous. This definition
clearly makes $y$ continuous on the open set $Q_0$, so it extends
continuously to the closure of $Q_0$. If we now extend $y$ to the whole of
$Q$ by making it zero on $Q\setminus\overline{Q_0}$ then we now have an
element of $C_\infty(Q)$ which is clearly the required supremum of the
set $A$.
\qed

\medskip

A characterization of Levi norms in terms of laterally increasing
sets was obtained  in [A1, Theorem~3$'$] or [A2, Theorem~2.5]. It is
the equivalence of (a) and (b) in our next Theorem~{2.3}.
However, for our further  work we  need  slightly more,
namely the equivalence of (a) and (c).

\medskip
\noindent{\bf Theorem 2.3.} {\sl For any Banach lattice $E$ the
following three properties are equivalent.

\item{\rm (a)}$E$ has a Levi norm.

\item{\rm (b)}Every norm bounded laterally increasing subset of
$E_+$ has a supremum.

\item{\rm (c)}If $A\subset E_+$ is a set of pairwise disjoint elements
          such that the set
$$
B=\big\{\sum_{a\in \sigma} a:\  \sigma \
                 {\rm is \ a\ finite\ subset\ of\ } A\big\}
$$
is norm bounded, then the set $A$ has a supremum (which will also be the
supremum of $B$).}

\smallskip
\noindent{\bf Proof.}
Implication (a)\implies(b)  is obvious. Since the proof in [A1] of the
implication   (b)\implies(a) is not easily available we, answering the request
of the referee,  present here a rather complete sketch of this proof.

Note first that (b) certainly implies that every order bounded set of
pairwise disjoint positive elements of $E$ has a supremum,
and so the theorem of  Veksler and Gejler cited above
tells us that $E$ must be Dedekind complete. This allows us to
embed $E$ as an order dense ideal in its universal completion
$\hat E = C_\infty (Q)$, where $Q$ is the Stonean space of $E$.

Let $(x_\alpha)$ be an increasing norm bounded net in $E_+$. We need to
consider the following two exclusive cases: either $(x_\alpha)$ is
order bounded in
$C_\infty (Q)$ or else
$(x_\alpha)$ is not bounded. In the former case there exists
$\ z=\sup_\alpha x_\alpha \in C_\infty (Q)$. Let $G_\alpha =
\overline{ \{q\in Q: 2x_\alpha(q) > z(q)\} }$. This is an open and closed set
for each $\alpha$, and clearly $G_{\alpha_1} \subset G_{\alpha_2}$ whenever
$\alpha_1 < \alpha_2 $. Therefore, the net
$(z\chi_{{}_{G_\alpha}})_{{}_\alpha}$ is laterally increasing. Since
$x_\alpha \uparrow z$ we have that $\bigcup G_\alpha$ is dense in $Q$, and this
implies that $z\chi_{{}_{G_\alpha}} \uparrow z$. It remains to notice that each
element
$z\chi_{{}_{G_\alpha}}$ belongs to $E$ since $z\chi_{{}_{G_\alpha}} \le
2x_\alpha$.  Consequently (b) implies that $z\in E$, that is, indeed, the net
$(x_\alpha)$  has its supremum in $E$.

Consider the second case when $(x_\alpha)$ is not order bounded in
$C_\infty (Q)$. Then there exists a nonempty open and closed subset
$Q_0$ of $Q$ such that \ $\sup_\alpha x_\alpha(q) = \infty$ for all $q\in D$,
where $D$ is a dense subset of $Q$. Take any $0\le z \in C_\infty(Q)$ with its
support in $Q_0$. Consider the net $(z\land x_\alpha)$.
This is an increasing norm bounded net in $E$, and clearly its supremum \
$\sup_\alpha z\land x_\alpha$ in $C_\infty(Q)$ exists and equals $z$, since \
$\sup_\alpha (z\land x_\alpha)(q) =z(q)$ for each $q\in D$.
By the previous part $z\in E$. In other words, we have proved that the
universally complete band $C_\infty(Q_0)$ is normable. This is clearly
impossible.

It is obvious that either of the conditions (a) or (b) implies (c),  by
considering the upward directed set $B$. We will prove that
(c)\implies(b), which will complete the proof of the equivalence of the
three statements.

Again notice that every  order bounded set of pairwise disjoint
positive elements of $E$ will certainly have a supremum by (c),
and so another application of the theorem of Veksler and Gejler
tells us that $E$ must be Dedekind complete. As before, we  assume
that $E$ is embedded into its universal completion $\hat E=C_\infty(Q)$.

   Let $A\subset E_+$ be laterally increasing and norm bounded, we
must show that $A$ has a supremum in $E$. Note that by Proposition 2.2,
$A$ has a supremum in $\hat E$, which we will denote by $y$. We claim
that $y\in E$. Without loss of generality we may assume that
$y(q)=1$ for each $q \in Q$, otherwise we will consider an appropriate
principal ideal generated in $C_\infty(Q)$ by the element $y$.
As above, we do not distinguish between members of the
ideal generated in $\hat E$ by $y$ and functions in  $C(Q)$.

The argument
used in the proof of Proposition 2.2 shows that on a dense subset $Q_0$ of
$Q$, we have  $1=y(s)=a(s)$ whenever $a\in A$ and $a(s)>0$. In particular,
this shows that each $a\in A$ is the characteristic function of some  open
and closed subset of $Q$ and that the union of these open and closed sets is
dense in $Q$. Consider now the collection of all open and closed subsets
$G$ of $Q$  each of which is   contained in some set
$\{s\in Q: a(s)=1\}$,
where $a$, of course,  depends on $G$.
Let ${\cal C}$ be a maximal disjoint collection
of such sets $G$. If $\bigcup \{G: G\in {\cal C}\}$
is not dense in $Q$ then there is $t\in Q_0$ which does not meet its
closure. But for some $a\in A$ we have $a(t)=1$. Adding the set
$\{s:a(s)=1\}\setminus\overline{\bigcup \{G: G\in {\cal C}\} }$ to
${\cal C}$  gives us a contradiction to the maximality of ${\cal C}$.
Thus the  family  $D=\{\chi_G:G\in {\cal C} \}$  is a subset of $C(Q)$
with supremum $y$.  Notice also that in actuality  each function $\chi_G$
from $D$  belongs to $E$, as $E$ is an ideal in $\hat E$ and $0\le
\chi_G\le a$ whenever $G\subseteq \{s:a(s)=1\}$.

If $\chi_{G_k}\in D$ for $k=1,2,\ldots,n$ then there are $a_k\in A$ with
$\chi_{G_k}\le a_k$ for $1\le k\le n$. As $A$ is upwards directed, there
is $b\in A$ with $a_k\le b$ for $1\le k\le n$. Hence $\chi_{G_k}\le b$
for each $k$. But the functions $\chi_{G_k}$ are disjoint, so we also
have $\sum_{k=1}^n \chi_{G_k}\le b$ and hence
$\|\sum_{k=1}^n \chi_{G_k}\|\le \|b\|$
so we may apply (c) to deduce that the family $D$ has a supremum,
$z$, in $E$. Clearly $z\le y$, so we may regard $z$ as an element
of $C(Q)$. But in $C(Q)$ it is clear that $z$ must be at least 1
on $\bigcup \{G: G\in {\cal C} \}$  which is dense in $Q$, so that $z\ge
y$.  That is, $z=y$, showing that $y\in E$ and, hence, that $A$ does
indeed have a supremum in $E$ (and not just in $\hat E$).
\qed

\medskip
For the sequentially Levi property we have the following analogue of the
previous theorem.

\medskip
\noindent{\bf Theorem 2.4.} {\sl For any
Banach lattice $E$ the following three properties are equivalent.

\item{\rm (a)}$E$ has a sequentially Levi norm.

\item{\rm (b)}Every norm bounded laterally increasing sequence in
          $E_+$ has a supremum.

\item{\rm (c)}Every disjoint positive sequence, for which the set of all
          possible finite sums is norm

          bounded, must have a supremum.}

\smallskip
\noindent The equivalence of (a) and (b) was established in [A1, Theorem~3]
or [A2, Theorem~{2.4}], while the equivalence of (c) and (b) is
obvious in view of the comment made after Definition~{2.1}.

We conclude this section by a simple observation.
{\sl If $X$ is a Dedekind $\sigma$-complete Banach lattice such that every
disjoint
family in $X$ is countable, then $X$ has a Levi norm if and only if $X$ has a
sequentially Levi norm}. Indeed, the disjoint families that need to be
considered in
Theorem 2.3 (c)  will be countable and then we can use Theorem 2.4, taking
into account
the fact that $X$ must be Dedekind complete by [V, Theorem VI.2.1].

The condition that every disjoint family in $X$ is countable is stronger than
the countable sup property, which also implies that $X$ is Dedekind complete.
However, the countable sup property alone is not enough to imply the equivalence
of Levi and sequentially Levi properties, as the Dedekind $\sigma$-complete
Banach
lattice, under the uniform norm, of all functions on [0,1] with at most
countable
support shows.

\bigskip

\noindent{\bf 3. Regularity of operators on L-spaces.}
As a first step in our proof of the converse of Theorem~{1.1} we need an
argument that involves only  separable domains,  so we isolate this as a
separate result.

\medskip
\noindent{\bf Theorem 3.1.}
{\sl The following conditions on a Dedekind $\sigma$-complete Banach
lattice $F$ are equivalent.

\item{\rm (a)}$F$ has a sequentially Levi norm.

\item{\rm (b)}For every separable L-space $E$ the equality ${\cal
L}(E,F)={\cal L}^r(E,F)$
holds.

\item{\rm (c)}For $E=L_1[0,2\pi]$ the equality ${\cal L}(E,F)={\cal
L}^r(E,F)$ holds.}

\smallskip
\noindent{\bf Proof.} (a)\implies(b).
Assume that $E$ is a separable L-space, $F$ has a sequentially Levi
norm and $U:E\to F$ is norm bounded. As in the usual proof
of  Theorem~{1.1} (see for example [V, Theorem~8.7.2]), for each
$0\le x\in E$, we consider the set
$$
A_x=\{\sum_{k=1}^n |Ux_k|: \ x_k\ge 0, \ \sum_{k=1}^n x_k=x,  \ n\in\N\}.
$$
This set is upward directed and norm bounded.
If  $A_x$  has a supremum,  then  by  the  classical Riesz-Kantorovich
formula  this supremum  is $|U|(x)$. In the familiar situation when $F$
has a Levi norm, this implies immediately that
$A_x$  has a supremum and the existence of $|U|$ follows.

In the present situation, however, when we have only a  sequentially Levi
norm, the existence of $\sup A_x$ is not obvious. To establish it we will
use an argument first utilized in  [W1,~Theorem~{5.2}]. To prove the
existence of
$\sup A_x$, it
will suffice to show that there is a {\it countable}  upward directed
subset $B_x$ of $F_+$ such that  $B_x$ is dense in $A_x$ and $B_x$ is
dominated by $A_x$, i.e., for each $b\in B_x$ there is $a\in A_x$ with
$b\le a$.  The latter
condition implies that $B_x$ is norm bounded, and so in view of the
sequentially Levi property this subset  $B_x$ will have a supremum, which
clearly will be also the supremum of $A_x$.

In order to find such a set $B_x$ it is sufficient to find a dense
countable subset, $C_x$, of $A_x$. Indeed,  the collection of all
finite suprema from that set $C_x$ will be clearly  upward directed,
countable and dense in $A_x$. It  will be  also dominated by $A_x$ as if
$y_1,y_2,\ldots,y_m\in  C_x\subseteq A_x$, then there is $z\in A_x$ with
$y_1,y_2,\ldots,y_m\le z$ (since $A_x$ is upward directed) and hence
$y_1\vee y_2\vee\cdots\vee y_m\le z$.

Finally, to prove the existence of  a countable dense subset of $A_x$, it
suffices (by standard arguments) to find a countable subset $D_x$ of $F$
with $A_x\subseteq\overline{D_x}$ (so we do not need this countable set to
be contained in $A_x$).
Since $E$ is separable, we can find a countable dense subset
$$E_x=\{z_k:k\in\N\}$$
of the order interval $[0,x]$. The set
$$D_x=\{\sum_{k=1}^n |Uz_{m_k}|: n\in\N\}$$
will certainly be countable, we show that its closure contains $A_x$
(and will in general be much larger).
Given a typical element $\sum_{k=1}^n
|Ux_k|$ of $A_x$ and $\epsilon>0$, for each $k$ we can find $z_{m_k}\in
E_x$ with $\|z_{m_k}-x_k\|<\epsilon/(n \|U\|)$. It follows that
$$
\big\||Ux_k|-|Uz_{m_k}|\big\|\le\|Uz_{m_k}-Ux_k\|
\le\|U\|\|z_{m_k}-x_k\|<\epsilon/n
$$
so that
$$\left\|\sum_{k=1}^n|Uz_{m_k}|-\sum_{k=1}^n|Ux_k|\right\|<\epsilon$$
showing that $A_x\subseteq\overline{D_x}$ as claimed.

Implication (b)\implies(c) is obvious. In order to prove that
(c)\implies(a), let us suppose, contrary to what we claim,
that $(e_n)$ is a disjoint positive sequence in $F$ such that

{(i)} $\left\|\sum_{n\in \sigma}e_n\right\|\le 1$ for all finite subsets
      $\sigma\subset \N$,

\noindent but with

{(ii)} $\{e_n:n\in\N\}$ not being bounded above

\noindent (if it were, then the set would
have a supremum as we are assuming that $F$ is Dedekind
$\sigma$-complete). We take $E=L_1[0,2\pi]$ and define
$$b_n(f)=\int_0^{2\pi} f(t) \cos(nt)\;dt$$
for each $f\in E$ and $n\in\N$. Clearly

{(iii)} $|b_n(f)|\le\|f\|_{{}_1}$,

\noindent and by the Riemann-Lebesgue theorem we know that

{(iv)} the sequence  $\big(b_n(f)\big)\in c_0$.

\noindent Define a linear operator $S:E\to F$ by $Sf=\sum_{n=1}^\infty
b_n(f) e_n$. It is routine to show, given (i) and (iv), that this series
converges in $F$, whilst the use of (i) and (iii) shows that $\|S\|\le 1$.

We claim that $S$ cannot be regular. If it were, let $T:E\to F$
be a positive operator
with $T\ge S$. We denote by $\bf 1$ the constantly one function on
$[0,2\pi]$, and let ${\bf 2}=2\cdot {\bf 1}$.

For each $n\in\N$ we have
$ {\bf 2}\ge \big({\bf 1}+\cos(nt)\big)\ge 0$,
and so
$$\eqalign{T{\bf 2}&\ge  T\big({\bf 1}+\cos(nt)\big)\cr
                    &\ge  S\big({\bf 1}+\cos(nt)\big)\cr
                    &= S{\bf 1}+ S\big(\cos(nt)\big)\cr
                    &= 0+ \pi\cdot e_n \ge 2e_n.\cr}$$
Thus, $T{\bf 1}\ge e_n$ for all $n\in\N$, which contradicts (ii).
\qed

\medskip
\noindent{\bf Remark 3.2.} The assumption that $F$ be Dedekind
$\sigma$-complete cannot be omitted from the hypotheses of Theorem 3.1.

Indeed it follows from [AG] that if $K$ is a Stonean space and $k_0\in K$
is an arbitrary non-isolated point then, for an arbitrary Banach lattice
$E$, every bounded operator from $E$ to $ C(\tilde K)$ is regular,
where the Hausdorff compact space $\tilde K$ is obtained from
the Hausdorff compact space $K\times\{1,2\}$ by identifying $(k_0,1)$ and
$(k_0,2)$. If we take $K=\beta(\N)$ and $k_0\in \beta(\N)\setminus\N$, then
$\N$ is a dense open $F_\sigma$-set in $K$ and clearly $k_0\in\overline{\N}$.
Therefore, the two sets
$\N\times\{1\}$ and $\N\times\{2\}$ are disjoint open $F_\sigma$-sets
in
$\tilde K$ and the intersection of their closures is non-empty since it
contains $(k_0,1)=(k_0,2)$. Hence $\tilde K$ is not even an
F-space, let alone a quasi-Stonean space, and consequently $C(\tilde K)$ is
not Dedekind $\sigma$-complete (in fact not even a Cantor space).
\qed
\medskip

So far we discussed the implications of the equality $(\star)$
on the properties of the target space $F$. There is one more
natural ``parameter" which has a very important impact on $F$,
namely  the order properties of the space ${\cal L}^r(E,F)$.
Following the lead given in [AG], we showed in [AW] that it is
possible to characterize  Dedekind $\sigma$-complete Banach
lattices $F$ by the fact that the space of regular operators from
any separable Banach lattice into $F$ forms a lattice.
To parallel that result in the setting of the present paper we need to
consider operators on separable L-spaces. To do so we will generalize
Theorem~3.10 in [AW], in which we considered operators defined on all
separable Banach lattices. We are taking this opportunity to also mention that
Theorem~{3.10}  of [AW] was rather carelessly worded (references to Banach
lattices
of operators should be replaced by references to vector lattices of operators;
the former were not removed from an earlier draft which was formulated
in a different setting), so we restate the result in its entirety,
together with   a new extra  equivalence.

\medskip
\noindent{\bf Theorem 3.3.} {\sl For a Banach lattice $F$ the following
are equivalent:
\item{\rm (a)} $F$ is Dedekind $\sigma$-complete.
\item{\rm (b)} For all separable Banach lattices $E$, ${\cal L}^r(E,F)$ is
a Dedekind
$\sigma$-complete vector lattice.
\item{\rm (c)} ${\cal L}^r(c,F)$ is a vector lattice, where $c$ is the
space of all
convergent sequences.
\item{\rm (d)} ${\cal L}^r(L_1[0,2\pi],F)$ is a vector lattice.}

\smallskip
\noindent{\bf Proof.} The equivalence of (a), (b) and (c) is proved
in Theorem 3.10 of [AW], whilst (d) is an obvious consequence of (b).

To complete the proof we will deduce (a) from (d). Notice first that if
${\cal L}^r(L_1[0,2\pi],F)$ is a vector lattice, then so is
${\cal L}^r(L_1[0,2\pi],J)$, where $J$ is any principal ideal in $F$. If we can
prove that each such $J$ is Dedekind $\sigma$-complete, then clearly, $F$
will also be. We may thus, using the Kakutani-Krein representation for
unital M-spaces, reduce the problem to that of showing that if $K$ is a
compact Hausdorff space and ${\cal L}^r(L_1[0,2\pi],C(K))$ is a vector lattice,
then $K$ is quasi-Stonean.

Let $U$ be an open $F_\sigma$-subset of $K$. By Proposition 2.1 of [W2]
we can find two disjoint sequences of functions in $C(K)$, $(u_n)$ and
$(v_n)$, vanishing on $K\setminus U$, such that $0\le u_n(k),v_n(k)\le 1$
for all $k\in K$ and $U=\bigcup_{n=1}^\infty \big(u_n^{-1}(1)\cup
v_n^{-1}(1)\big)$. If $k_0\in K\setminus \overline{U}$, use Urysohn's
lemma to find $w\in C(K)$ with $0\le w(k)\le 1$ for all $k\in K$,
$w_{|U}\equiv1$ and $w(k_0)=0$. With the same definition of
$b_n$ as in Theorem 3.1, define two linear operators
$S,T:L_1[0,2\pi]\to C(K)$ by
$$
\eqalign{Sf&=\sum_{n=1}^\infty b_n(f) u_n\cr
         Tf&=Sf+b_0(f)w,\cr}
$$
 where $b_0(f):=\int_0^{2\pi} f(t)\,dt$.
The convergence of the series defining $S$ may be proved in a similar
manner to the  proof in Theorem~{3.1}. We notice for later use
that the operator $S$ is independent of our choice of the point $k_0$ and of
the function $w$.

It is clear that if $f\ge 0$ then $(T-S)f=b_0(f)w\ge 0$, so that
$T\ge S$. Next we claim that the operator $T$ is positive.
To this end first note that $b_0+b_m$ is a positive functional
for each $m\in \N$. Indeed
$b_0(f)+b_m(f) = \int_0^{2\pi} f(t) (1+\cos(mt))\,dt \ge 0$,
whenever $f\ge 0$. We have
$$
Tf=\sum_{n=1}^\infty b_n(f) u_n +b_0(f)w = b_m(f)u_m + b_0(f)w +
\sum_{n\neq m} b_n(f) u_n.
$$
Since $w\ge u_n\ge 0$  and  the functions $(u_n)$ are pairwise disjoint,
the positivity of the functionals $b_0+b_m$ implies
that $Tf\ge 0$ on the closure of the set $\bigcup_n \{k\in K:u_n(k)>0\}$.
On the complement of the above set (i.e., on the interior of the set where all
$u_n$ vanish), $Tf$ is simply $b_0(f)w$ which is certainly non-negative, so
that $Tf\ge 0$. Thus we have established that $T\ge S,0$ which, among other
things, tells us that $S$ is a regular operator. By condition (d), $S^+$ exists
in ${\cal L}^r(L_1[0,2\pi],C(K))$.

Using again the inequality ${\bf 2}\ge {\bf 1}+\cos(nt)\ge0$ we see that
$$
S^+{\bf 2}\ge S^+\big({\bf 1}+\cos(nt)\big)
\ge S\big({\bf 1}+\cos(nt)\big)= \pi\cdot u_n,
$$
whence $\pi^{-1}S^+{\bf 2}$ is an upper bound for the sequence $(u_n)$. We must
also have $S^+\le T$, so that $0\le S^+{\bf 2}\le T{\bf 2}=2w$ and in particular
$(S^+{\bf 2})(k_0)=0$. Since our operator $S$ was independent of
the choice of the point $k_0$, the previous argument is applicable to
any point of $K\setminus \overline{U}$. Consequently,
we must have $S^+{\bf 2}$ identically zero on
$K\setminus\overline{U}$.

If we replace the sequence $(u_n)$ by $(v_n)$  to define an operator
$S_1f=\sum_{n=1}^\infty b_n(f) v_n$ and  repeat the whole of the
proof so far, then we will conclude that $\pi^{-1}S_1^+{\bf 2}$
vanishes on $K\setminus\overline{U}$ and is an upper  bound for $(v_n)$.

Now $\pi^{-1}S^+{\bf 2}\vee S_1^+{\bf 2}$ also vanishes on
$K\setminus\overline{U}$ and is an upper bound for both sequences
$(u_n)$ and $(v_n)$. It follows that
$\pi^{-1} S^+{\bf 2}\vee S_1^+{\bf 2}$  is at least 1 on $U$, whilst
it vanishes on $K\setminus\overline{U}$.
Since $S^+{\bf 2}\vee S_1^+{\bf 2}$
is continuous, this certainly  forces
$\overline{U}$ to be open, and hence $K$
is indeed quasi-Stonean.
\qed

\medskip
Putting together Theorems 3.1 and 3.3 we now have:

\medskip
\noindent{\bf Corollary 3.4.} {\sl The following conditions on a
Banach lattice $F$ are equivalent.
\item{\rm (a)}$F$ has a sequentially Levi norm.
\item{\rm (b)}For every separable L-space $E$,  the space
${\cal L}(E,F)$ is a vector lattice.} \medskip
We turn now to the general case, when we allow the domain to be
any L-space instead of  considering only separable L-spaces. Our
results are entirely analogous to those above and, in fact, depend
on them.

\medskip
\noindent{\bf Theorem 3.5.} {\sl The following conditions on a Dedekind
complete Banach lattice $F$ are equivalent.

\item{\rm (a)}$F$ has a Levi norm.

\item{\rm (b)}For every L-space $E$, the equality ${\cal L}(E,F)={\cal
L}^r(E,F)$ holds.}

\smallskip
\noindent{\bf Proof.} As mentioned earlier the
implication (a)\implies(b) is true by  Theorem~{1.1}.
We need only prove that (b)\implies(a). Suppose
that $A\subset F_+$ is a disjoint set
such that for all finite sets $\sigma\subset A$,
$\left\|\sum_{a\in \sigma} a\right\|\le K$.
By Theorem 2.3, it is enough to show that $A$ has
a supremum in $F$. Let $\alpha$ denote the cardinality of $A$. We will work
with operators whose domain is the space $L_1(\mu)$, where $\mu$ is the
product of $\alpha$ copies of the probability measure on $\{0,1\}$ which
assigns measure ${1\over 2}$ to both $\{0\}$ and $\{1\}$. It is well known
that integrable functions on $\{0,1\}^\alpha$
depend on only countably many
variables. Let $\phi_i$ denote the function in $L_\infty(\mu)$ which
depends only on the $i$'th variable and takes the value $1$ if this variable
is $0$, and takes the value $-1$ if this variable is $1$. If $f\in L_1(\mu)$
and  $f$ does not depend on the $i$'th variable, then
$\int f \phi_i\;d\mu=0$.
 We denote  by $\Phi$
the collection of all these functions $\phi_i$ for $i\in \alpha$.
Note that we may also write  $A=\{a_i:i\in \alpha\}$ by
indexing the members of $A$ by $\alpha$. We thus certainly have

   {(i)} For each $f\in L_1(\mu)$, $\int f \phi\;d\mu=0$ for all but
countably many $\phi\in\Phi$.

  {(ii)} $\|\phi\|_\infty=1$ for all $\phi\in\Phi$.

\noindent For notational convenience, we will write $\phi(f)$ in place of
$\int f\phi\;d\mu$ and regard each $\phi$ as an element of $L_1(\mu)^*$.
Note that

  {(iii)} for each $\phi\in \Phi$ there is $f\in L_1(\mu)$ with
$0\le f\le {\bf 1}$ and $|\phi(f)|\ge {1\over 2}$.

\noindent We refer the reader to [HS, \S22] for the requisite details
concerning infinite product measures and integration.

In view of  Theorem 3.1 we already know that $F$ has a sequentially Levi
norm. This implies, by a theorem of Amemiya [Am], that there
is a constant $C>0$ such that $0\le y_n\uparrow y\Rightarrow
\|y\|\le C \lim \|y_n\|$. In particular, for each countable subset
$B\subset A$ we have that  its supremum $\bigvee B$ exists
in  $F$ and that $\left\|\bigvee B\right\|\le C$.
Define a linear operator $S:L_1(\mu)\to F$ by
$$Sf=\sum_{i\in\alpha} \phi_i(f) a_i.$$

In order to see that this series is order convergent, recall that (i)
guarantees that there is a countable subset $\beta\subset\alpha$ such
that $\phi_i(f)=0$ for all $i\in\alpha\setminus\beta$. The collection of
all finite sums $\sum_{i\in \sigma}\phi_i(f) a_i$ is norm bounded
 in view of the inequality
$$
{\left|\sum_{i\in \sigma} \phi_i(f) a_i\right|\le
\sum_{i\in \sigma}\|f\|_1 a_i =\|f\|_1\sum_{i\in \sigma} a_i}
$$
which implies
that $\left\| \sum_{i\in \sigma} \phi_i(f) a_i\right\|\le K\|f\|_1$. The
disjointness of the  elements
$a_i$ guarantees that we also
have $\left\|\sum_{i\in \sigma}\phi_i(f)^\pm a_i\right\|\le K\|f\|_1$, so it
follows from the sequentially Levi property of $F$ that the series
$\sum_{i\in\beta}\phi_i(f)^\pm a_i$ are both order convergent and hence
$\sum_{i\in\beta}\phi_i(f) a_i$ is also order convergent. This implies
that $Sf$ is indeed well-defined.

Notice that Amemiya's result shows us that $\left\|\sum_{i\in
\beta}\phi_i(f)^\pm a_i\right\|\le C K \|f\|_1$, so that
$\|Sf\|=\left\|\sum_{i\in \beta}\phi_i(f) a_i\right\|\le 2C K \|f\|_1$,
so that $\|S\|\le 2 C K$ and, in particular, $S$ is norm bounded.
By (b) there is $T:L_1(\mu)\to F$ with $T\ge S,-S$.
We know that for each $i\in\alpha$  there is $f_i\in L_1(\mu)$ with
$0\le f_i\le {\bf 1}$ and such that $|\phi_i(f_i)|\ge {1\over 2}$. Thus
$$
T{\bf 1}\ge Tf_i\ge|Sf_i|\ge |\phi_i(f_i)| a_i\ge {1\over2}a_i,
$$
so that $T{\bf 2}$ is an upper bound for $A$. As we are assuming that $F$
is Dedekind complete, this implies that the supremum of $A$ exists.
\qed

\medskip
The example given in  Remark 3.2 shows that we cannot omit
the assumption of Dedekind completeness from the statement of
Theorem~{3.5}. The L-space $E$ produced in the proof above is a
nonseparable $L_1(\mu)$-space with a finite measure $\mu$. It is
interesting to notice that one cannot avoid using a somewhat extravagant
measure space to get the desired contradiction.
For example, the classical L-spaces $\ell_1(\Gamma)$ will be definitely not
enough, since  for any $\Gamma$ and any Banach lattice $F$ each continuous
operator from  $\ell_1(\Gamma) $ into  $F$ is regular.

Similarly to what was done in Theorem~{3.3}, our next result
characterizes the Dedekind completeness of $F$ in terms
of order properties of the space ${\cal L}^r(E,F)$.

\medskip
\noindent{\bf Theorem 3.6.} {\sl For a Banach lattice $F$ the following
are equivalent:
\item{\rm (a)}$F$ is Dedekind complete.
\item{\rm (b)}For all Banach lattices $E$, the space ${\cal L}^r(E,F)$ is a
Dedekind
complete  vector lattice.
\item{\rm (c)}For all L-spaces $E$, the space ${\cal L}^r(E,F)$ is a vector
lattice.}

\smallskip
\noindent{\bf Proof.} Again, it is only (c)\implies(a) that we need
to prove.
As in the proof of Theorem~{3.3}, it suffices to consider the case that
$F=C(K)$. By that theorem, we already know that $C(K)$ is Dedekind
$\sigma$-complete, and so $K$ is quasi-Stonean.

Let $U$ be an arbitrary open subset of $K$. We need to prove that
its closure $\overline{U}$ is open. There obviously exists
a maximal disjoint collection of closed and open subsets
$D_i \  (i\in I)$ of $U$, and so $\cup D_i$ is dense in $U$.
Fix an arbitrary point   $k_0\in K\setminus\overline{U}$, and  find
a function $w\in C(K)$ which lies between {\bf 0} and {\bf 1}, is
one on $U$ and  is zero at the point $k_0$.

Exactly as  in the proof of Theorem~{3.5}, we
construct now a measure $\mu$ on $\{0,1\}^\alpha$,
where $\alpha$ is the cardinality  of $I$.  Let $\Phi$ have
the same meaning as in that proof and let $\psi$ be the linear
functional $f\mapsto \int f\;d\mu$, so that $\psi-\phi_i\ge0$ for all
$i\in I$.
Define $S,T: L_1(\mu)\to C(K)$
by
$$
\eqalign{Sf&=\sum_{i\in I}\phi_i(f) \chi_{{\!}_{D_i}}\cr
Tf&=Sf+\psi(f)w.\cr}
$$
This time convergence of the series follows from the Dedekind
$\sigma$-completeness of $F$, using the facts that for each
$f\in L_1(\mu)$ only countably
many terms are non-zero and  that $\psi\ge \phi_i$ for all $i\in
\alpha$, and showing that the finite partial sums all lie between
$\pm\psi(f) w$. It should be pointed out that the operator $S$ above
is independent of the function $w$ and of the point $k_0$. As in the
proof of Theorem 3.1, we can show that $T\ge S,0$,
so that $S^+$ exists  by (c). For each $i\in I$ there
is $f_i\in L_1(\mu)$ with $0\le f_i\le {\bf 1}$ and $\phi_i(f_i)\ge
{1\over2}$. It follows  as before
that $S^+{\bf 2}\ge 2S^+f_i\ge \chi_{{\!}_{D_i}}$.
Noting that $\phi_i({\bf 1})=0$ we also have
$0\le S^+{\bf 2}\le T{\bf 2}=2w$. In particular,
$S^+{\bf 2}(k_0)=0$. As $k_0$  was an arbitrary point
of $K\setminus\overline{U}$, this shows that $S^+{\bf 2}$ must vanish
on $K\setminus\overline{U}$. On the other hand,
the continuous function $S^+{\bf 2}$  dominates each
$\chi_{{\!}_{D_i}}$, so is at least one on a dense subset of
$U$. It follows that
$\overline{U}$ is indeed open and the proof is complete.
\qed

\medskip
\noindent{\bf Corollary 3.7.} {\sl The following conditions on a Banach
lattice $F$ are equivalent.
\item{\rm (a)} $F$ has a Levi norm.
\item{\rm (b)} For every L-space $E$,  the space ${\cal L}(E,F)$ is a
vector lattice.}
\medskip

The results obtained in this section imply the
following observation. If a Dedekind complete Banach lattice $F$
does not have a Levi-norm then, by Theorem~{3.5}, for an appropriate
L-space $E$ the space ${\cal L}(E,F)$ is bigger than ${\cal L}^r(E,F)$, though
the latter is a Dedekind complete vector lattice.

\bigskip

\noindent{\bf 4. Regularity of operators with arbitrary domain.}
There is only one known case of a Dedekind complete Banach lattice $F$
such  that all
continuous operators, with any Banach lattice as domain
 and $F$ as range,
are regular---namely when $F$ has a strong order unit. This result
dates back to [K]. As we shall see next, it is in fact the only case.
\medskip
\noindent{\bf Theorem 4.1.} {\sl The following conditions on a Dedekind
complete Banach
lattice $F$ are equivalent.
\item{\rm (a)} $F$ has a strong order unit.
\item{\rm (b)} For every Banach lattice $E$ the equality ${\cal
L}(E,F)={\cal L}^r(E,F)$
holds.}

\smallskip
\noindent{\bf Proof.} As we said (a)\implies(b) is due to Kantorovich. We
need only prove that (b)\implies(a). First notice that by Theorem 3.5,
$F$ has a Levi norm.
Second  note  that (b) is ``more" than enough to imply that $F$ is
isomorphic to an M-space.
Indeed, either of the following two conditions is weaker than (b) and
implies that $F$ is isomorphic to an M-space.

1) There exists $p\in (1,\infty]$ such that each continuous operator
$T: L_p[0,1]\to F$ is regular.

2) There exists a Banach lattice $E$ containing uniformly the subspaces
$\ell_n^\infty$ and such that  each continuous operator $T: E\to F$ is
regular.

The sufficiency of the first condition was established by
Cartwright and Lotz [CL] (see also [MN, Theorem~3.2.1], or
[A3, Theorem~8.6]), whilst for sufficiency of the second and several
other conditions see [A3, \S2]. The  conjunction of the two  properties of
$F$ obtained so far immediately implies (a). Indeed, the collection of all
finite suprema from the unit ball of $F$ is upward directed and norm
bounded, so has a supremum. That supremum is a strong order unit for $F$.
\qed

\medskip
Again,
as  Remark~{3.2} shows,
the hypothesis of Dedekind completeness may not be omitted.
Combining Theorem 3.5 and Theorem 4.1, we obtain:

\medskip
\noindent{\bf Corollary 4.2.} {\sl The following conditions on a Banach
lattice $F$ are equivalent.
\item{\rm (a)} $F$ is Dedekind complete and has a strong order unit.
\item{\rm (b)} For every Banach lattice $E$, the space ${\cal L}(E,F)$ is a
vector
lattice.}
\medskip

We conclude by noticing that in all our results starting with Corollary 3.4
the uniform operator norm and the regular norm on the space
${\cal L}^r$ are equivalent. The isometric version of these results,
describing  when these two norms are equal, will be discussed in
a forthcoming paper by the authors and Z. L. Chen.

\vfill
\eject

\bigskip
\centerline {\bf References}
\bigskip

[A1] Y. A. Abramovich,
            Some theorems about normed lattices,
           {\it Vestnik Leningr. Univ. Mat. Meh. Astronom.\/  No.
            \bf 13}$\,$(1971), 5--11.

[A2] Y. A. Abramovich,
           {\it On some properties of norm in normed lattices
           and  their maximal normed extensions}, Thesis, Leningrad State
           Univ. 1972.

[A3]  Y. A. Abramovich,
            When each continuous operator is regular,
           {\it Funct. Analysis, Optimization and Math. Economics, Oxford
            Univ. Press}, 1990, 133--140.

[AG] Y. A. Abramovich (Ju. A. Abramovi\v c), V. A. Gejler,
            On a question of Fremlin concerning order bounded and regular
            operators,  {\it Coll. Math.,\/ } {\bf 46}$\,$(1982), 15--17.

[AW]  Y. A. Abramovich and A. W. Wickstead,
       The regularity of order bounded operators into C(K), II,
       {\it Quart. J. Math. Oxford \/\bf 44}$\,$(1993), 257--270.

[AB] C. D. Aliprantis and  O. Burkinshaw,
               {\it Positive Operators\/}, Pure and Applied
               Mathematics Series \#$\,$119, Academic Press,
               New York \& London, 1985.

[Am] I. Amemiya,
       A generalization of Riesz-Fischer's theorem,
       {\it J. Math.Soc. Japan \/ \bf 5}$\,$(1953), 353--354.

[CL] D.I. Cartwright and H.P. Lotz,
             Some characterizations of AM- and AL-spaces, {\it Math. Z.\/
             \bf 142}$\,$(1975) 97--103.

[F] D. Fremlin, {\it Topological Riesz Spaces and Measure Theory},
                Cambridge Univ. Press, London, 1974.

[HS] E. Hewitt and K. Stromberg,
              {\it Real and Abstract Analysis},  Springer-Verlag,
               Berlin Heidelberg New York, 1969.

[K] L. V. Kantorovich,
          Concerning the general theory of operators in semi-ordered spaces,
          {\it Dokl. Acad. Nauk USSR,\/ \bf  1}$\,$(1936), 271--274.

[KV] L. V. Kantorovich, B. Z. Vulikh,
           Sur la repr\'esentation des op\'erations lin\'eares,
           {\it Compos. Math. \/\bf 5}$\,$(1937), 119--165.

[MN] P. Meyer-Nieberg, {\it Banach Lattices}, Springer-Verlag,
            Berlin Heidelberg New York, 1991.

[MN1] P. Meyer-Nieberg, Characterisierung einiger topologischer und
         ordnungstheoretischer Eigenschaften von Banachverbanden mit Hilfe
         disjunkter Folgen, Arch. Math. (Basel) 24 (1973), 640-647.

[N] H. Nakano, {\it Modulared Semi-Ordered Linear Spaces}, Maruzen,
        Tokyo, 1950.

[S] J. Synnatzschke,
           On the conjugate of a regular operator and some
           applications to the questions of complete continuity of regular
           operators,  {\it Vestnik Leningr. Univ. Mat. Meh. Astronom.\/
           n.\bf 1}$\,$(1972), 60--69.

[VG] A. I. Veksler and V. A. Gejler,
     Order and disjoint completeness of linear partially ordered spaces,
     {\it Sibirsk. Math. J.\/}  {\bf 13}$\,$(1972), 43--51.

[V] B. Z. Vulikh,
    {\it Introduction to the theory of partially ordered
     spaces,\/ } Walters-Noordoff Sci. Publication, Groningen,  1967.

[W1] A. W. Wickstead,
     The regularity of order bounded operators into $C(K)$,
     {\it Quart. J. Math. Oxford}$\,$(2) {\bf 41}$\,$(1990), 359--368.

[W2] A. W. Wickstead,
           Spaces of operators with the Riesz Separation Property,
           {\it Indag. Math.\/ \bf 6}$\,$(1995), 235--245.

[Zaa] A. C. Zaanen, {\it Riesz Spaces} II, North Holland, Amsterdam, 1983.

\vskip 1.4cm

\settabs  8 \columns

\+{\bf  Y. A. Abramovich}     &&&&  {\bf A. W. Wickstead} \cr
\+ Department of Mathematics  &&&& Department of Pure Mathematics \cr
\+ IUPUI 		                   &&&& The Queen's University of Belfast \cr
\+ Indianapolis, IN 46202     &&&& Belfast BT7 1NN  \cr
\+ USA                        &&&& Northern Ireland  \cr

\bye